\documentclass[12pt]{iopart}
\usepackage{amssymb}

\begin{document}

\note[Semi-symmetric spacetimes]{Investigation of all Ricci
semi-symmetric and all conformally semi-symmetric spacetimes}

\author{Jan E. {\AA}man}

\address{Department of Physics, Stockholm University, SE-106 91 \ Stockholm,
Sweden}
\eads{\mailto{ja@fysik.su.se}}

\begin{abstract}

We find all Ricci semi-symmetric as well as all conformally
semi-symmetric spacetimes. Neither of these properties implies the
other. We verify that only conformally flat spacetimes can be Ricci
semi-symmetric without being conformally semi-symmetric and show that
only vacuum spacetimes and spacetimes with just a $\Lambda$-term can
be Ricci semi-symmetric without being conformally semi-symmetric.

\end{abstract}

\pacs{04.20.Cv, 02.40.Ky, 04.20.Jb}

\section{Introduction}

Semi-symmetric spaces where introduced by Cartan \cite{Cartan1946} and
are characterized by the curvature condition
\begin{equation} \label{SemiSymmetric}
    \nabla_{[a}\nabla_{b]}R_{cdef} = 0,
\end{equation}
where $R_{cdef}$ denotes the Riemann tensor and round and square
brackets enclosing indices indicate symmetrization and
antisymmetrization, respectively.

A semi-Riemannian manifold is said to be {\em conformally
semi-symmetric} if the Weyl tensor $C_{abcd}$ satisfies
\begin{equation} \label{ConformallySemiSymmetric}
    \nabla_{[a}\nabla_{b]}C_{cdef} = 0;
\end{equation}
and {\em Ricci semi-symmetric} if the Ricci tensor satisfies
\begin{equation} \label{RicciSemiSymmetric}
    \nabla_{[a}\nabla_{b]}R_{cd} = 0.
\end{equation}

The geometrical properties of semi-symmetric spaces where discussed
in \cite{Sz,HV,Senovilla2006,Senovilla2008,DDVV}.

In this paper we find all semi-symmetric, Ricci semi-symmetric and all
conformally semi-symmetric spacetimes. Eriksson and Senovilla
\cite{Eriksson2010} found all such non-conformally flat spacetimes. 

It is an advantage to instead of tensors rather to use the Weyl spinor
$\Psi_{ABCD}$, the curvature scalar spinor $\Lambda$, and the spinor
$\Phi_{ABA'B'}$ for the tracefree part of the Ricci tensor.  We also
use the spinor ${\rm X}_{ABCD} = \Psi_{ABCD} +\Lambda \left(
\varepsilon_{AC}\varepsilon_{BD} +\varepsilon_{AD}\varepsilon_{BC}
\right)$. The spinor commutator $\Box_{AB}$ operating on spinors with
one index is $\Box_{AB}\, \kappa_C = -{\rm X}_{ABC}{}^E\kappa_E$ and
$\Box_{AB}\, \tau_{C'} = -\Phi_{ABC'}{}^{E'}\tau_{E'}$ \cite{Penrose1984}.

In spinors (\ref{ConformallySemiSymmetric}) is equivalent to
$\Box_{AB} \Psi_{CDEF} = 0$ and $\Box_{A'B'} \Psi_{CDEF} = 0$, or
\begin{eqnarray}
     {\rm X}_{AB(C}{}^G \Psi_{DEF)G} &= 0, \label{WeylCondition1} \\
      \Phi_{A'B'(C}{}^G \Psi_{DEF)G} &= 0. \label{WeylCondition2}
\end{eqnarray}

\section{Conditions for semi-symmetry}

Calculation of the components of (\ref{WeylCondition1}) in terms of
$\Psi_{ABCD}$ and $\Lambda$ shows show that it not will have 15
independent components but rather just 5 components and
(\ref{WeylCondition1}) can be replaced by its contraction over $BC$ or
\begin{equation}
\Psi^{GH}{}_{(AD}\Psi_{EF)GH} -2\Lambda \Psi_{ADEF} = 0.
\label{PsiCondition1}
\end{equation}

We observe that the spinor commutator $\Box_{AB}$ yields 0 if
operating on a scalar as the curvature scalar $\Lambda$, therefore
only its effect on $\Phi_{ABA'B'}$, i.e. on the tracefree part of the
Ricci tensor has to be considered.  The second condition for Ricci
semi-symmetry (\ref{WeylCondition2}) corresponds to spinor equations
\begin{eqnarray}
   \Box_{AB} \Phi_{CDC'D'} &=
        -2 {\rm X}_{AB(C}{}^E \Phi_{D)EC'D'}
        -2 \Phi_{AB(C'}{}^{E'} \Phi_{D')E'CD}, \label{RicciCondition1} \\
   \Box_{A'B'} \Phi_{CDC'D'} &=
        -2 \overline{\rm X}_{A'B'(C'}{}^{E'} \Phi_{D')E'CD}
        -2 \Phi_{A'B'(C}{}^E \Phi_{D)EC'D'}. \label{RicciCondition2}
\end{eqnarray}
The condition (\ref{RicciCondition2}) is however just the complex
conjugate of (\ref{RicciCondition1}) and will not give any additional
conditions for Ricci semi-symmetry.

Equation (\ref{RicciCondition1}) has three groups of two symmetric
indices but will in fact not have 27 independent components as
$\Box^{AB} \Phi_{ABC'D'} = 0$, it can be replaced by the fully
symmetric spinors
\begin{eqnarray}
 \Box_{(AB} \Phi_{CD)C'D'} = - 2 \Psi_{(ABC}{}^E\Phi_{D)EC'D'}
  \label{PhiCondition1} \\
 \Box_{(A}{}^F \Phi_{C)FC'D'} = 4 \Lambda\Phi_{ACC'D'}
   -\Psi^{EF}{}_{AC}\Phi_{EFC'D'}
   -2 \Phi^E{}_A{}^{F'}{}_{(C'}\Phi_{D')F'CE} \label{PhiCondition2}
\end{eqnarray}
with 15 and 9 components respectively.

The various spacetimes are for simplicity studied in a frame where
first $\Psi_{ABCD}$ has been brought to a standard form depending on
its Petrov type. Thereafter $\Phi_{ABA'B'}$ is brought to standard
form depending on its Segre type \cite{Quartic,Exact}.

All of the above formulas have been implemented in CLASSI
\cite{Classi}.

\section{Main Results}

The equations for Ricci semi-symmetric spacetimes
(\ref{PhiCondition1}) and (\ref{PhiCondition2}) and Petrov Types
\textbf{I}, \textbf{II} or \textbf{III} in standard frame requires
that $\Phi_{ABA'B'}=0$, while spacetimes with a $\Lambda$-term only or
vacuum ($\Lambda=0$) are Ricci semi-symmetric.

For Petrov type \textbf{D} spacetimes (with only $\Psi_2$ nonzero) the
conditions both for conformal and for Ricci semi-symmetry require all
components of $\Phi_{ABA'B'}$ to be zero except for $\Phi_{11'}$. The
conditions for conformal semi-symmetry reduces to $2 \Lambda +\Psi_2 =
0$ and for Ricci semi-symmetry to $\Phi_{11'} (2 \Lambda +\Psi_2) =
0$. Spacetimes of Segre type A1[(11)(1,1)] are therefore both Ricci
and conformally semi-symmetric if they conform to
$\Lambda=-\frac{1}{2}\Psi_2$ with arbitrary $\Phi_{11'}$.  Spacetimes
with $\Lambda$-term only are semi-symmetric if once again
$\Lambda=-\frac{1}{2}\Psi_2$ but only Ricci semi-symmetric for other
relations between $\Psi_2$ and $\Lambda$.

For Petrov type \textbf{N} spacetimes (with only $\Psi_4$ nonzero) the
conditions for Ricci semi-symmetry as well as for conformal
sem-symmetry require all components of $\Phi_{ABA'B'}$ to be zero
except for $\Phi_{22'}$. The conditions here reduces to $\Lambda
\Psi_4 = 0$ for conformal semi-symmetry and $\Lambda \Phi_{11'} = 0$
for Ricci semi-symmetry. They are therefore both Ricci and conformally
semi-symmetric for arbitrary $\Psi_4$ and $\Phi_{22'}$ as long as
$\Lambda=0$. Once again all spacetimes with $\Lambda$-term only are
Ricci semi-symmetric.

\subsection{Conformally flat spacetimes}

This is the case not treated by Eriksson and Senovilla \cite{Eriksson2010}.

All conformally flat (Petrov type \textbf{0}) spacetimes are
conformally semi-symmetric.

Calculations with CLASSI shows that all conformally flat spacetimes of
Segre type A1[(11)(1,1)] (only $\Phi_{11'}$ nonzero) with $\Lambda=0$
are also Ricci semi-symmetric.

For conformally flat perfect fluids (Segre type A1[(111),1)],
$\frac{1}{2}\Phi_{00'}=\Phi_{11'}=\frac{1}{2}\Phi_{22'}$) as well as for
tachyons (Segre type A1[1(11,1)],
$-\frac{1}{2}\Phi_{00'}=\Phi_{11'}=-\frac{1}{2}\Phi_{22'}$) the
condition for Ricci semi-symmetry reduces to $\Lambda = \Phi_{11'}$.
So perfect fluids with
$\Lambda=\frac{1}{2}\Phi_{00'}=\Phi_{11'}=\frac{1}{2}\Phi_{22'}$ as
well as tachyons with
$\Lambda=-\frac{1}{2}\Phi_{00'}=\Phi_{11'}=-\frac{1}{2}\Phi_{22'}$ are
also Ricci semi-symmetric.

All spacetimes of Segre type A3[(11,2)] with $\Lambda=0$ (only
$\Phi_{22'}$ nonzero) are also Ricci semi-symmetric.

All conformally flat spacetimes with $\Lambda$-term only, including
flat spacetimes, are Ricci semi-symmetric.

For all other Segre types the conditions for Ricci semi-symmetry are
not satisfied.

\section{Summary}

All conformally flat spacetimes are conformally semi-symmetric, all
spacetimes with $\Lambda$-term only are Ricci semi-symmetric.

Semi-symmetric spacetimes are flat spacetimes, $\Lambda$-term and
Segre A1[(11)(1,1)] spacetimes of Petrov type \textbf{0} and of type
\textbf{D} with $\Lambda=-\frac{1}{2}\Psi_2$, Segre type A3[(11,2)]
with $\Lambda=0$ of Petrov types \textbf{N} and \textbf{0}, as well as
conformally flat perfect fluids and tachyons with
$\Lambda=\Phi_{11'}$. For a table see appendix A.

\section*{Acknowledgments}

I thank Brian Edgar (1945-2010) for pointing out ref
\cite{Eriksson2010} and introducing me to semi-symmetric spacetimes. I
thank Ingemar Bengtsson and Jos\'e Senovilla for comments on the
manuscript.

\newpage

\section*{Appendix A}
\begin{table}[h]

\begin{tabular}{|l|c|c|c|c|c|c|}
\hline
Segre type \ \ \ \ \ \ \ \ \ Petrov type & \textbf{I} & \textbf{II} &
  \textbf{III} & \textbf{D} & \textbf{N} & \textbf{0} \\
\hline
$\Lambda$-term A1[(111,1)] or vacuum & Ric s-s & Ric s-s & Ric s-s &
  Ric s-s & Ric s-s & semi-sym \\
$\Lambda$-term, $\Lambda=-\frac{1}{2}\Psi_2$ & Ric s-s & Ric s-s &
  $\nexists$ & semi-sym & $\nexists$ & $\nexists$ \\
A1[(11)(1,1)], $\Lambda=-\frac{1}{2}\Psi_2$ & - & - & $\nexists$ & semi-sym &
  $\nexists$ & $\nexists$ \\
A1[(11)(1,1)] & - & - & - & see above & - & semi-sym \\
A3[(11,2)], $\Lambda=0$ & - & - & - & - & semi-sym & semi-sym \\
A1[(111),1] perfect fluid, &  &  &  &  &  & \\
  \ \ $\Lambda=\frac{1}{2}\Phi_{00'}=\Phi_{11'}=\frac{1}{2}\Phi_{22'}$
  & - & - & - & - & - & semi-sym \\
A1[1(11,1)] tachyon fluid, &  &  &  &  &  & \\
  \ \ $\Lambda=-\frac{1}{2}\Phi_{00'}=\Phi_{11'}=-\frac{1}{2}\Phi_{22'}$ &
  - & - & - & - & - & semi-sym \\
All other Ricci tensors & - & - & - & - & - & conf s-s \\
\hline
\end{tabular}

\caption{Relations between Petrov type, Segre type and conformal
semi-symmetry (conf s-s), Ricci semi-symmetry (Ric s-s) and
semi-symmetry. A hyphen (-) indicates neither conformal nor Ricci
semi-symmetry.}

\end{table}

\end{document}